\documentclass{amsart}

\usepackage{graphicx}
\usepackage{amsfonts}
\usepackage{amsmath}
\usepackage[ruled,vlined,norelsize]{algorithm2e}
\usepackage{colortbl}
\usepackage{commath}

\newtheorem{theorem}{Theorem}[section]

\theoremstyle{definition}
\newtheorem{definition}[theorem]{Definition}
\newtheorem{example}[theorem]{Example}

\theoremstyle{remark}

\SetAlFnt{\small}
\SetAlCapFnt{\small}
\SetAlCapNameFnt{\small}
\SetAlCapHSkip{0pt}
\IncMargin{-\parindent}
\numberwithin{equation}{section}

\begin{document}

\title{Computing the Braid Monodromy of Completely Reducible $n$-gonal Curves}

\author{MEHMET E. AKTAS}
\address{Department of Mathematics, Florida State University, Tallahassee, Florida 32306}
\email{maktas@math.fsu.edu}

\author{ESRA AKBAS}
\address{Department of Computer Science, Florida State University, Tallahassee, Florida 32306}
\email{akbas@cs.fsu.edu}





\begin{abstract}
Braid monodromy is an important tool for computing invariants of curves and surfaces. In this paper, the \emph{rectangular braid diagram (RBD)} method is proposed to compute the braid monodromy of a completely reducible $n$-gonal curve, i.e. the curves in the form $(y-y_1(x))...(y-y_n(x))=0$ where $n\in \mathbb{Z}^{+}$ and $y_i\in \mathbb{C}[x]$. Also, an algorithm is presented to compute the Alexander polynomial of these curve complements using Burau representations of braid groups. Examples for each computation are provided.
\end{abstract}

\maketitle

\section{Introduction}

Braid monodromy is an important tool for computing invariants of curves and surfaces. It is used in the study of the topology of plane curve complements to compute the fundamental group \cite{zariski1929problem},\cite{van1933fundamental} and explicitly defined in \cite{chisini1933suggestiva}. It is also used for covering spaces \cite{chisini1952courbes} and Alexander polynomial \cite{Libgober1989}. In \cite{moishezon1991braid}, it is found that the braid monodromy defines not just the fundamental group but also the homotopy type. A. Libgober also relates braid monodromy and homotopy type of curve complements \cite{libgober1986homotopy}. In \cite{kulikov2000braid}, braid monodromy is also used to analyze connected components of the moduli space of surfaces of general type. The authors also had important results for the relationship between topology and braid monodromy in the particular case of nodal-cuspidal curves.  In \cite{carmona2003monodromia}, the author generalized Kulikov-Teicher results. There are also works relating topology and braid monodromy by Artal-Carmona-Cogolludo \textit{et al.} \cite{bartolo2003braid,artal2007effective,artal2001sextics,artal2005topology}. Kharlamov-Kulikov studied on braid monodromy factorization \cite{kharlamov2003braid}. Degtyarev used the Grothendieck's \textit{dessin d'enfants} to compute the braid monodromy of trigonal curves and the plane curves with degree $d$ that have a singular point of multiplicity $d-3$ which are birationally equivalent to trigonal curves \cite{degtyarev2012topology}. 

Although braid monodromy has a simple construction and is a widely used invariant, explicit computational algorithms have only been developed for some special cases. In \cite{moishezon1981stable}, authors had an algorithm to compute the braid monodromy of real singular curves that have a cusp point. Oka computed the braid monodromy for some special cases \cite{oka1992symmetric,oka1996two,oka2002fundamental}. Bessis had more general algorithm used for exceptional braid groups \cite{bessis2004explicit}. This algorithm resulted the VKCURVE package in GAP.  For the complexified real arrangements, the braid monodromy is computed in \cite{hironaka1993abelian,salvetti1987topology,cordovil1995braid}. Arvola computed the braid monodromy of an arbitrary complex arrangement \cite{arvola1992fundamental}. In \cite{kaplan2004braid}, authors improved an algorithm for an almost real curve. Artal-Carmona-Cogolludo had also algorithms for the case of real curves \cite{artal2005topology,artal2007effective}. In \cite{deconinck2001computing}, authors computed the monodromy of an irreducible algebraic curve which led to the \textit{algcurves} package in Maple. Amor{\'o}s computed the braid monodromy of branching divisor $R \in \mathbb{C}^2$ using the vertical projection based on numerical integration of a related differential equation \cite{amoros2012monodromias}.

The Alexander polynomial is an important topological invariant for curve complements. For example, in Zariski's famous example, the sextic with 6 cusps where the cusps are on conic has the Alexander polynomial $t^2+t+1$ whereas the sextic with 6 cusps where the cusps are not on conic has $1$. Although the Alexander polynomial is very useful, there are few algorithms and implementations in literature to compute the Alexander polynomials of curve complements. In \cite{eko}, the Alexander polynomial of $M$ in the case when $M$ is a 3-manifold fibered over a circle is computed using monodromies. In \cite{csimcsek2004computer}, the authors wrote a computer program to calculate the Alexander polynomial from braid monodromy of a given knot. In both papers, they used Fox calculus and assumed that they know the braid monodromies already. 

An $n$-gonal curve is a compact algebraic curve with a linear pencil of degree $n$. It can also be defined as a plane algebraic curve given by $F(x,y)=0$ where deg$_yF=n$. Moreover, $n$-gonal curves are closely related with the plane curves with degree $d$ that have a singular point of multiplicity $d-n$ where $d>n$. Because, if one blows-up that singular point, the proper transform of the plane curve is an $n$-gonal curve in a rational ruled surface with an exceptional section. 

In this paper, we propose a new method, named as the \textbf{R}ectangular \textbf{B}raid \textbf{D}iagram (RBD), to compute the braid monodromy of a completely reducible $n$-gonal curve. In this method, we construct the rectangular braid diagram of a given completely reducible $n$-gonal curve. Using the diagram, we get the loops around its singular fibers and compute the braid monodromy on each loop using an adaptive step size method. Experiments show that RBD method is more efficient and also comprehensive than the VKCURVE \cite{bessis2004explicit}. In addition to computing the braid monodromy, we present an algorithm to compute the Alexander polynomial of $n$-gonal curves based on Libgober's theorem in \cite{Libgober1989}. Hence, using RBD and Libgober's theorem, the Alexander polynomial of completely reducible $n$-gonal curves can be computed by just using their defining polynomials. We furthermore give some examples for each computation.

The paper is structured as follows: In Section 2, we give preliminaries and definitions. In Section 3, we introduce the RBD algorithm for braid monodromy computations and go over some examples. In Section 4, we present our algorithm that computes the Alexander polynomial and continue with some examples. In Section 5, we compare RBD method with VKCURVE. Finally, conclusions are drawn in Section 6.

\section{Preliminaries}

\subsection{The $n$-gonal curves}

Let $\Sigma=\Sigma_k, k\geq 0$, be a Hirzebruch surface i.e. a rational geometrically ruled surface with an exceptional section $E=E_k$ of self-intersection $-k$. Denote the ruling by $p:\Sigma\rightarrow B$ where the base $B$ is a genus 0 curve, i.e. $B \simeq \mathbb{P}^1$. For a point $b$ in $B$, denote the fiber $p^{-1}(b)$ by $F_b$. 

\begin{definition}
An \emph{$n$-gonal curve} is a reduced curve $C\in \Sigma$ not containing $E$ or a fiber of $\Sigma$ as a component such that the restriction $p: C \rightarrow B$ is a map of degree $n$ i.e. each fiber intersects with $C$ at $n$ points counting with multiplicities. In other words, they are the plane algebraic curves given by $F(x,y)=0$ where deg$_y F=n$. 
\end{definition}

A \emph{singular fiber} of an $n$-gonal curve $C\subset \Sigma$ is a fiber $F$ of $\Sigma$ intersecting $C + E$ in less than $n+1$ distinct points. Hence, $F$ is singular either it passes through $C\cap E$, or $C$ is tangent to $F$ or $C$ has a singular point in $F$.  

\begin{definition}
An $n$-gonal curve $C$ is \textit{completely reducible} if it is defined by $(y-y_1(x))...(y-y_n(x))=0$ where $y_i \in \mathbb{C}[x]$ for all $i\in \{1,...,n\}$ and $n \in \mathbb{Z}^+$.
\end{definition}

\subsection{The Braid Group $\mathbb{B}_n$}

Let $\mathfrak{F}_n=\langle \alpha_1, ... , \alpha_n \rangle$ be the free group on $n$ generators. The \emph{braid group} $\mathbb{B}_n$ can be defined as the group of automorphisms $\beta: \mathfrak{F}_n\rightarrow \mathfrak{F}_n$ with the following properties:
\begin{itemize}
  \item each generator $\alpha_i$ is taken to a conjugate of a generator;
  \item the element $\rho:= \alpha_1...\alpha_n$ remains fixed.
\end{itemize}

In \cite{artin1947theory}, Artin showed that $\mathbb{B}_n=\langle \sigma_1,...,\sigma_{n-1}|\sigma_i\sigma_{i+1}\sigma_i=\sigma_{i+1}\sigma_i\sigma_{i+1}, [\sigma_i,\sigma_j]=1$ if $|i-j|>1\rangle$ and the action of $\mathbb{B}_n$ on $\mathfrak{F}_n$ as follows:

\begin{center}
$\sigma_i(\alpha_j)=
    \begin{cases}
      \alpha_{j} & j \neq i,i+1 \\
      \alpha_i\alpha_{i+1}{\alpha_i}^{-1} & j=i \\
      \alpha_i & j=i+1
\end{cases}
$
\end{center}
One of the oldest and most well-known representation of the braid group is discovered by Burau \cite{burau1935zopfgruppen}. The Burau representation $\rho$  of $\mathbb{B}_n$ can be described as the $\mathbb{Z}[t,t^{-1}]$ representation which maps the generators of $\mathbb{B}_n$ as follows:\\
\begin{flushleft}
$\sigma_1 \rightarrow$ $\left(
\begin{tabular}{c c|c}
  $-t$ & 1 & 0 \\
  0 & 1 & 0 \\
  \hline
  0 & 0 & $I_{n-3}$
 \end{tabular}\right)$
\\
\vspace{1mm}
$\sigma_i \rightarrow$ $\left(
\begin{tabular}{c|c c c|c}
  $I_{i-2}$ & 0 & 0 & 0 & 0\\
  \hline
  0 & 1 & 0 & 0 & 0\\
  0 & $t$ & $-t$ & 1 & 0\\
  0 & 0 & 0 & 1 & 0\\
  \hline
  0 & 0 & 0 & 0 & $I_{n-i-2}$
 \end{tabular}\right)$, $2\leq i \leq n-2$,
\\
\vspace{1mm}
$
\sigma_{n-1} \rightarrow$ $\left(
\begin{tabular}{c | c c}
  $I_{n-3}$ & 0 & 0 \\
  \hline
  0 & 1 & 0 \\
  0 & $t$ & $-t$
 \end{tabular}\right)$
\end{flushleft}

This representation will be used in Chapter 4 to compute the Alexander polynomial of the $n$-gonal curve complements.

\subsection{The Braid Monodromy}
Let $C\subset \Sigma$ be an $n$-gonal curve and $p:\Sigma \rightarrow B$ be the ruling. Let $F_1, F_2, ..., F_r$ be the singular fibers of $C$ and $E$ be the exceptional section. Denote $b_i=p(F_i)$, the image under the ruling of the corresponding singular fiber $F_i$. 

Take a monodromy domain $\Lambda\subset B$ where $\Lambda=B\setminus D^{\circ}$ and $D$ is an embedded closed disk. WLOG, we can assume that the interior of $\Lambda$ contains all $b_i$'s. Pick a point $b\in B$ such that the fiber $p^{-1}(b)=F$ is a nonsingular fiber. Let $F^{\sharp} =F\setminus (C\cup E.)$. Clearly, $F^{\sharp}$ is equal to $F\setminus E$ with $n$ punctures. Let $\Lambda^{\sharp} = \Lambda \setminus \{ b_1, b_2,..., b_r\}$.

We know that $\pi_1(F^{\sharp})$ is the free group on $n$ generators i.e. $\pi_1(F^{\sharp})=\mathfrak{F}_n= \langle \alpha_1, ..., \alpha_n \rangle$ where $\alpha_i$ is the loop which covers $i$-th intersection of the fiber $F^{\sharp}$ and the $n$-gonal curve $C$. We should note that many different generator systems and their relations with braids exist for $\pi_1(F^{\sharp})$. A precise description about this will be done in Section \ref{sec:3.1}. $\pi_1(B^{\sharp})$ is the free group on $r$ generators i.e. $\pi_1(B^{\sharp})=\mathfrak{F}_{r-1}=\langle \gamma_1,...,\gamma_{r-1} \rangle$ where $\gamma_j$ is the loop around $p_j$.

The restriction $p:p^{-1}(\Lambda \setminus (C\cup F))\rightarrow \Lambda^{\sharp}$ is a locally trivial fibration.	For each $j = 1,...,r$, dragging the fiber $F$ along $\gamma_j$ and keeping the base point results to the \textit{braid monodromy} homomorphism
\begin{center}
$\mathfrak{m}(\gamma_j):\pi_1(\Lambda^{\sharp},b)\rightarrow \text{Aut }\pi_1(F^{\sharp}).$ 
\end{center}
Since along any loop $\gamma \in \pi_1(\Lambda^{\sharp},b)$, the braid monodromy $\mathfrak{m}(\gamma)$ takes a generator $\alpha \in \pi_1(F^{\sharp})$ to a conjugate of a generator while keeping the product of all generators is fixed, $\mathfrak{m}(\gamma)$ is in the braid group $\mathbb{B}_n$. 

Each braid monodromy $\mathfrak{m}(\gamma_j)$ called \textit{the local braid monodromy} of the singular fiber $F_j$. The set of all local braid monodromies $\{\mathfrak{m}(\gamma_1),...,\mathfrak{m}(\gamma_r)\}$ is called the \textit{global braid monodromy} the $n$-gonal curve $C$. 
\section{Computing the Braid Monodromy}

In this section, we compute the global braid monodromy of the completely reducible $n$-gonal curve complements. Let $C\subset \Sigma \rightarrow B$ be a completely reducible $n$-gonal curve. The very first step to compute the braid monodromy is to find the singular fibers of the complement. These fibers are passing through either singular points of one irreducible component or intersections of these components. In our case, any irreducible component has no singular point by itself since ${\partial (y-y_i)}/{\partial y} \neq 0$. Hence we only need to find the intersection points of these $n$ components $\{y-y_i\}, i \in \{1,...,n\}$.

After finding all singular fibers, we choose a base point $b\in B$ that should not be from one of the singular points and then find loops starting from $b$ and inclosing each singular point. Next, we start to walk on each loop and find the change of the position of each strand, which comes from each irreducible component, and compute the corresponding braid element in $\mathbb{B}_n$ for each position change and finally get all braid monodromies for each loop (what we mean by ``position change" will be explained in coming section).

There are mainly two parts of this algorithm. The first part is to find the loops for each singular point. To find them, we construct the \emph{rectangular braid diagram}.

\subsection{The Rectangular Braid Diagram\label{sec:3.1}}

Let $V\subset \Sigma$ be the set of all singular points and $W=p(V) \subset B$, i.e. the image of $V$ under the ruling map $p: \Sigma \rightarrow B$. We will compute the monodromy at infinity $\mathfrak{m}_{\infty}$ separately at the end. Hence, if $\infty \in W$, take $W'=W \setminus \lbrace \infty \rbrace$, else, set $W'=W$. 

Since $W'\in \mathbb{C}$, one can talk about the real and imaginary parts of the points in $W'$. Let Re$(W')$ and Im$(W')$ be the set of all real and imaginary parts of each point in $W'$, respectively. There may be repeated values in these sets. For example, it is possible to have two conjugated singular points where they have the same real parts. First, we remove these repeated values in these sets. Then, we sort all points in Re$(W')$ and Im$(W')$. Then, we find the mid-points $r_i$ and $c_j$ of each consecutive two points in Re$(W')$ and Im$(W')$ respectively and divide the plane by $x=r_i$ and $y=c_j$ into rectangular parts. Moreover, we add the lines $x=\text{max}\{\text{Re}(W')\}+1, x=\text{min}\{\text{Re}(W')\}-1$, and $y=\text{max}\{\text{Im}(W')\}+1, y=\text{min}\{\text{Im}(W')\}-1$. As a result, we get the rectangular braid diagram. Note that there is at most one singular point in each rectangular region.

Now, we need to find the loop for each singular point. First, we choose a base point $b$ on the intersection of the right-most vertical line and $x$-axis, say $b=(x_0,0)$. We find the set of all rectangular regions that have a singular point in it, call it $Rec$. Let $b' \in W'$ and take the rectangle $\varsigma' \in Rec$ that incloses $b'$. To create the corresponding loop $\gamma'$, we start to move from the base point $b=(x_0,0)$ to the point $b_1=(x_0,y_{\text{top}})$ where $y_{\text{top}}$ is the imaginary value of the top-right point of $\varsigma$. Then we move from $b_1$ to $b_{\text{top}}=(x_{\text{top}},y_{\text{top}})$ where $b_{\text{top}}$ is the top-right point of $\varsigma$. After that, we continue to walk on $\varsigma$ in the counterclockwise direction and come back to $b_{\text{top}}$. We follow the inverse of our previous move and return to the base point $b$. As a consequence, we construct the loop $\gamma'$.
Furthermore, since $B\simeq \mathbb{P}^1$, i.e it is isomorphic to Riemann sphere, to compute the monodromy at infinity $\mathfrak{m}_{\infty}$, it is enough to take the loop that incloses $W'$. To do this, take the smallest rectangle that incloses all the rectangles in $Rec$, $\gamma_\infty$, and $\mathfrak{m}_{\infty}$ is equal to the monodromy around $\gamma_\infty$. Hence, we get all the loops based on $b$ and incloses each point in $W$.

Note that since each irreducible component has no singular point and the coefficient of $y^n$ in $C$ is equal to $1$ for all $x$, the analytic continuation holds in this construction. In other words, for every point $p$ on any loop, the corresponding fiber $F_p$ intersects with $C$ exactly at $n$ points. 

\begin{example}
Let $C$ be defined by $(y-x^2)(y-x-1)(y+1)=0$. Then its singular fibers are $p_1=-2$, $p_2=\frac{1+\sqrt{5}}{2}\approx 1.618$, $p_3=-i$, $p_4=i$, $p_5=\frac{1-\sqrt{5}}{2}\approx -0.618$ (our algorithm uses the approximations of algebraic numbers as it happens here) and its corresponding braid diagram is given in Figure 1.
\end{example}

\begin{figure}[ht!]
\centering
\includegraphics[width=10cm]{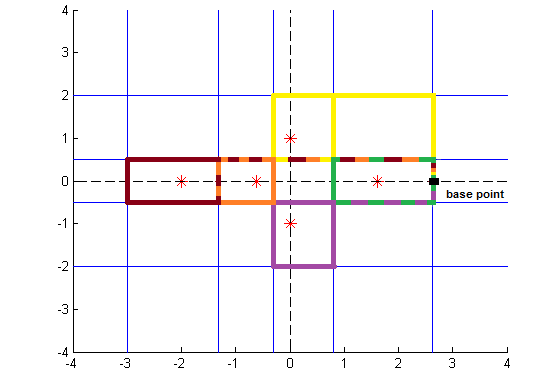}
\caption{The rectangular braid diagram of the curve $(y-x^2)(y-x-1)(y+1)=0$}
\end{figure}

\subsection{Computing the braid monodromy on a loop}

The second part of the RBD algorithm is computing the braid monodromy on each loop. To do this, we move on each loop with $\epsilon$ steps and in each step, we sort the real values of each strand and check whether the order changes or not (if the real values of two strands coincide, we assume the order does not change). If the order changes in a step, we find out at which strands it happens. This gives the braid element. We take the exponent of the braid element positive if the strands moves counterclockwise and take negative otherwise. To check this orientation, we look at the imaginary values of the strands.

For example, if the order of the real values of $k^{th}$ and $(k+1)^{st}$ strands changes, it results either the braid element $\sigma_k$ or  $\sigma_k^{-1}$. For the orientation, if the imaginary value of $k^{th}$ strand is smaller than the imaginary value of the $(k+1)^{st}$ strand, then there is a counterclockwise rotation which results $\sigma_k$. If not, it gives $\sigma_k^{-1}$. We repeat this step along each loop.

\subsubsection{Choosing the step size} 

The step size $\epsilon$ has a big importance in getting the correct results and also in efficiency of the algorithm. In some cases, the step size $\epsilon$ may not be small enough that it may miss some braid elements and we get wrong answers. On the other hand, the step size $\epsilon$ may be too small that causes unnecessary computations and it reduces the efficiency of algorithm. Even these two cases may happen on the same curve. Therefore, we need to choose $\epsilon$ for each step independently. Here is the algorithm to choose it:

Let $C:(y-y_1(x))...(y-y_n(x))=0$ be a completely reducible $n$-gonal curve, $\gamma$ is a loop in its rectangular braid diagram and $b$ be its base point. First, start with a value for the step size, let $\epsilon=e_0$. We need to investigate whether we miss braid elements on the interval $(b,b+\epsilon)$. To do this, we find upper bounds on how far each component $y_i$ for $i\in \{1,...,n\}$ moves on $(b,b+\epsilon)$. Let $y_i'=\sum_{j=1}^k a_jx^j$ where $a_j\in \mathbb{C}$ for $j\in \{1,...,k\}$. We have 
\begin{center}
$\left|y_i(b+\epsilon)-y_i(b)\right|=\left|\int_b^{b+\epsilon} y_i' dx\right|=\left|\int_b^{b+\epsilon} \sum_{j=0}^k a_jx^j dx \right |=\sum_{j=0}^k\left|  \int_b^{b+\epsilon}a_jx^j dx \right |\leq \sum_{j=0}^k  \int_b^{b+\epsilon}\left|a_j\right | \left| x^j\right| dx \leq \sum_{j=0}^k  \left|a_j \right| \text{max}\{\norm{b},\norm{b+\epsilon}\}^j\left|\epsilon\right|.$
\end{center}
Set $y_{i,\epsilon}=\sum_{i=1}^k  \left|b_i \right| \text{max}\{\norm{b},\norm{b+\epsilon}\}^i$ for all $i \in \{1,...,n\}$. We then check whether $y_{i,\epsilon}+y_{j,\epsilon}<\left | y_i(b)-y_j(b)\right |$ is true for all $1\leq i \leq j \leq n$. If they are all true, we guarantee that there is at most one braid element in each step. If not, set $\epsilon=e_0/2$ and do the some process till $y_{i,\epsilon}+y_{j,\epsilon}<\left | y_i(b)-y_j(b)\right |$ is true for all $1\leq i \leq j \leq n$. This is the last step of our algorithm. Algorithm \ref{alg} outlines the RBD algorithm.

\begin{algorithm}[ht!]\label{alg}
\DontPrintSemicolon
\KwData{$y_1(x),...,y_n(x)$ where $y_i\in \mathbb{C}[x]$ and the curve is $C:(y-y_1(x))...(y-y_n(x))=0$}
\KwResult{The braid monodromies of the complement of $C$}
\Begin{
$V \longleftarrow$ {roots of $(y_1-y_2),(y_1-y_3),...,(y_1-y_n),(y_2-y_3),...,(y_{n-1}-y_n)$}\;
Find the rectangular braid diagram\;
\For{$v\in V$}{
    Find the loop of $v$ in the braid diagram
}
$L \longleftarrow$ {set of loops of singular points}\;

\For{$l \in L$}{
$l_{s} \longleftarrow$ {set of sides of $l$}\;
\For{$ \iota\in l_s$}{
$b_{sta} \longleftarrow$ {starting point of $\iota$}\;
$b_{end} \longleftarrow$ {ending point of $\iota$}\;
$e_0 \longleftarrow$ {starting $\epsilon$ value}\;

$B \longleftarrow sorted \{y_1(b_{sta}),...,y_n(b_{sta})\}$\;

\While{$b_{sta}\neq b_{end}$}{
$\epsilon \longleftarrow$ {adjusted step size}\;

$B' \longleftarrow sorted\{y_1(b_{sta}),...,y_n(b_{sta})\}$\;
\If{$B\neq B'$}{
Find where $B(i)\neq B'(i), \forall i\in \{1,...,n\}$ and the corresponding braid
}
$b_{sta}=b_{sta}+\epsilon$
}
}

}
}
\caption{Braid Monodromy\label{BM}}
\end{algorithm}

\begin{example} 
Let $C$ be a trigonal curve defined by $(y-x)(y+x)(y-1)=0$. This curve has three singular fibers such as $x=-1,x=0,x=1$ and $x=\infty$. The corresponding braid monodromies are

\begin{center}
$m_{-1}={\sigma_2}^{-1}{\sigma_1}^{-1}{\sigma_2}^{2}{\sigma_1}{\sigma_2}$,

$m_{0}={\sigma_2}^{-1}{\sigma_1}^{2}{\sigma_2}$,

$m_{1}={\sigma_2}^{2}$,

$m_\infty=(\sigma_2\sigma_1\sigma_2)^2$.

\end{center}
\end{example}

\begin{example}
In \cite{amram2007fundamental}, the authors compute the local braid monodromy of the curve $(y+2x)(y+x^2)(y-x^2)$ at $(0,0)$. They state that two points (corresponding to $y=1$ and $y=-1$) do two counterclockwise full-twists and the third point (corresponds to $y=-2$) does counterclockwise full-twist around them, that is the braid monodromy is $m={\sigma_1}{\sigma_2}^{2}{\sigma_1}{\sigma_2}^{4}$. When we use our algorithm, we found the monodromy as 
$${\sigma_2}{\sigma_1}{\sigma_2}{\sigma_1}{\sigma_1}{\sigma_2}{\sigma_1}{\sigma_2}.$$ 
If we use the braid relation ${\sigma_2}{\sigma_1}{\sigma_2}={\sigma_1}{\sigma_2}{\sigma_1}$ four times, we get the same monodromy.
\end{example}

Notice that in the previous example, the authors could only compute the local braid monodromy but we can compute the global monodromies using our algorithm.

\section{Computing the Alexander polynomials}
In this section, we compute the Alexander polynomial of a completely reducible $n$-gonal curve using its global monodromy. Our main reference for this section is \cite{Libgober1989}.

\subsection{Algorithm and Implementation}

Let $C$ be an $n$-gonal curve and $\{p_1,...,p_r\}$ be images of the singular fibers of $C$ under the ruling map. Let $\rho$ be a $d$ dimensional linear representation of the braid group $\mathbb{B}_n$ over the ring $A$ of Laurent polynomials $\mathbb{Q}[t,t^{-1}]$. Define $P(C,\rho)$ as the greatest common divisor of the order $d$ minors in the $N\times d$ matrix of the map $\bigoplus (\rho(\mathfrak{m}(\gamma_i))-$Id$)$, where $N=rd$ , $\gamma_i$ is the loop enclosing $p_i$ and $\mathfrak{m}(\gamma_i)$ is the braid monodromy of the loop $\gamma_i$. We call this matrix the \emph{Libgober matrix}. It takes $(A^d)^{N}$ to $(A^d)$. Now, we can state the following important theorem that we will use in our algorithm.

\begin{theorem}\em \cite{Libgober1989}
If $\rho$ is the reduced Burau representation then $P(C,\rho)$ is equal to the Alexander polynomial $\bigtriangleup_C(t)$ of $C$ multiplied by $(1+t+...+t^{n-1})$.
\end{theorem}

This invariant is defined for plane algebraic curves that have generic projections and an $n$-gonal curve may not have a generic projection. However, we can still use this theorem to compute the Alexander polynomial because when $\rho$ is the reduced Burau representation, $\bigoplus (\rho(\mathfrak{m}(\gamma_i))-$Id$)$ is same with the matrix of Fox derivative $\bigoplus (\partial\mathfrak{m}(\gamma_i)g_k/\partial g_l-$Id$)$ where $g_k, g_l$ are generators of the fundamental group of the curve complement and Fox derivatives can be used in Alexander polynomial computations for any plane curves. 

In the previous section, we computed the global braid monodromy. Since the reduced Burau representation of $\mathbb{B}_n$ is $n-1$ dimensional, we need to find the $(n-1) \times N$ Libgober matrix of the braid monodromies and compute the set $M$ of $(n-1)\times (n-1)$ minors of this matrix. Then we find the $gcd$ of all elements in $M$ and divide it by $1+t+...+t^{n-1}$ to get the corresponding Alexander polynomial.

\subsection{Results}

First, we state one definition and one theorem that are important to interpret our results.

\begin{definition} 
The Alexander polynomial of a curve $C$ is defined to be trivial if $\bigtriangleup(t)=(t-1)^{r-1}$ where $r$ is the number of irreducible components of $C$.
\end{definition}

\begin{theorem}\em	\label{theorem}
Assume that $C_1$ and $C_2$ intersect transversely and let $C=C_1\cup C_2$. Then the Alexander polynomial $\bigtriangleup(t)$ of $C$ is given by $(t-1)^{r-1}$ where $r$ is the number of irreducible components of $C$, i.e., it is trivial \cite{Oka2005}. 
\end{theorem}

Now, we state some examples that we have found out by using our algorithm. We show how we get the results in the very first example only since we use the same steps for the others.

As the first example, we take a trigonal curve where each component intersects transversely:   

\begin{example}\label{example::1} 
The Alexander polynomial of the curve $(y-x^2)(y-x-1)(y-1)=0$, whose real picture is given in Figure \ref{fig:2}, is $(t-1)^2$.
\begin{figure}[ht!]
\centering
\includegraphics[width=6.5cm]{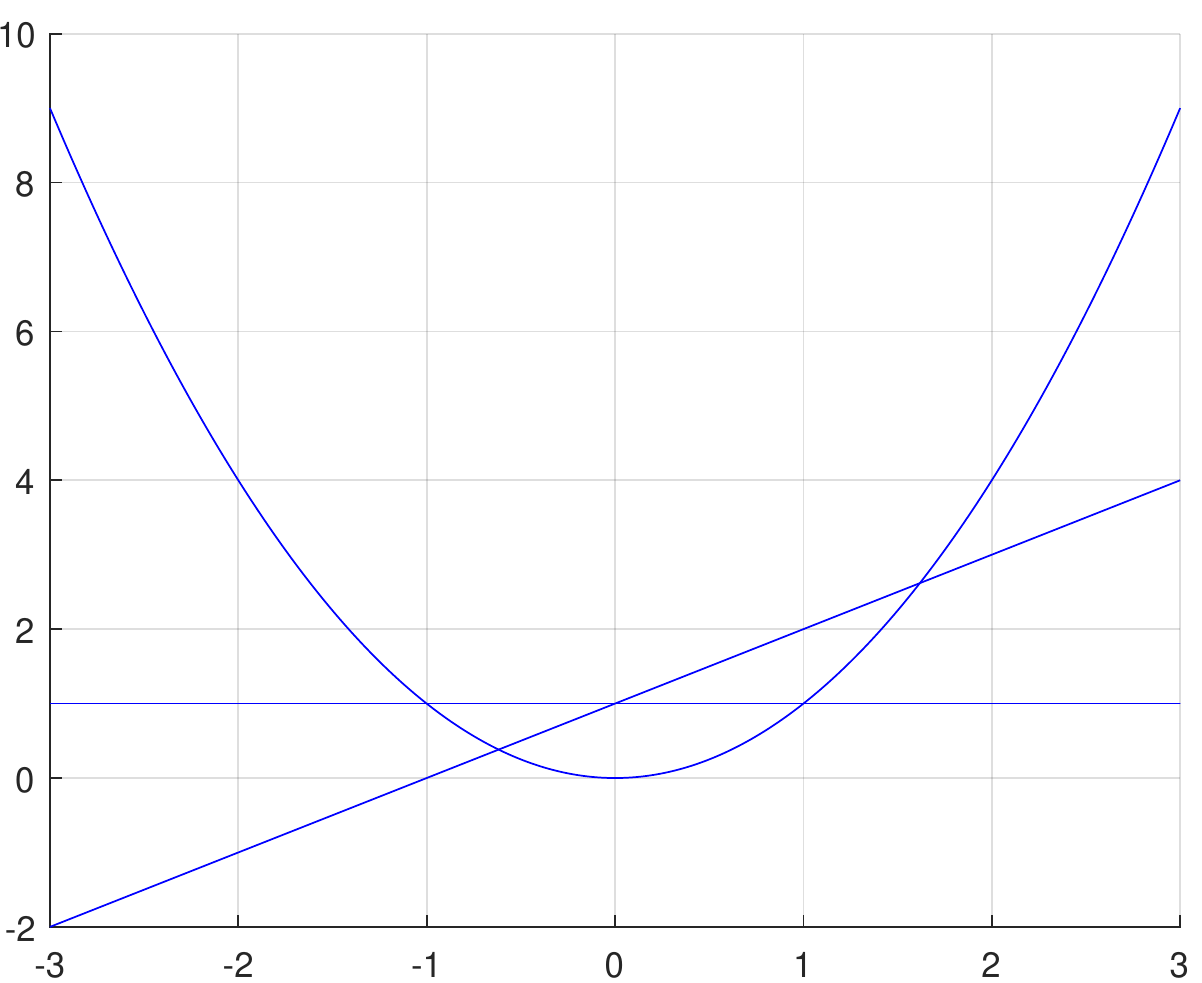}
\caption{The real picture of $(y-x^2)(y-x-1)(y-1)=0$}
\label{fig:2}
\end{figure}
\end{example}

\begin{proof}
There are totally six singular fibers of this curve at $W=\{0, -1, 1, \frac{-\sqrt{5} + 1}{2}, \frac{\sqrt{5} + 1}{2},\infty \}$ and, from the rectangular braid diagram technique, the braid monodromies are $\{\sigma_1^2,\sigma_2\sigma_1^2\sigma_2^{-1},\sigma_2\sigma_1^2\sigma_2^{-1},\sigma_2^2,\sigma_2^2,\sigma_2\sigma_2(\sigma_2\sigma_1\sigma_2)^2\sigma_1\sigma_2 \}$ respectively. 

Then, we compute the Libgober matrix using the braid monodromies:
\begin{center}
$\left(
\begin{tabular}{c c}
$t^2-1$ & $-t+1$ \\
0 & 0 \\
$t - 1$ &   $t^2 - t$\\
$t - 1$ &   $t^2 - t$\\
$t - 1$  & $t^2 - t$\\
$t - 1$  & $t^2 - t$\\
0  &       0\\
$-t^2 + t$&   $t^2 - 1$\\
0   &      0\\
$ -t^2 + t$ &   $t^2 - 1$\\
$t^4 - 1$ & $t^5 - t^4$\\
$0$ &  $t^6 - 1$\\
\end{tabular}\right).$

\end{center}

As the final step, find the $gcd$ of all order 2 minors in the Libgober matrix and divide it by the polynomial $t^2+t+1$ which gives the Alexander polynomial

\begin{center}
$\bigtriangleup(t)=(t-1)^2.$
\end{center}
\end{proof}

Note that in Example \ref{example::1}, the fundamental group also can easily be computed using its braid monodromy. In fact, it is abelian and this implies the triviality of the Alexander polynomial. However, in general, computing the fundamental group is not an easy task. That is why the Alexander polynomial is defined as an indeterminate invariant.

Furthermore, we realized that in some examples, Alexander polynomial is also trivial even if there are some non-transversal intersections, see the following example.

\begin{example} 
The Alexander polynomial of the curve $(y-x^2)(y-2x)(y+2x)(y)=0$, whose real picture is given in Figure \ref{fig:3}, is $(t-1)^3$.
\begin{figure}[ht!]
\centering
\includegraphics[width=6.5cm]{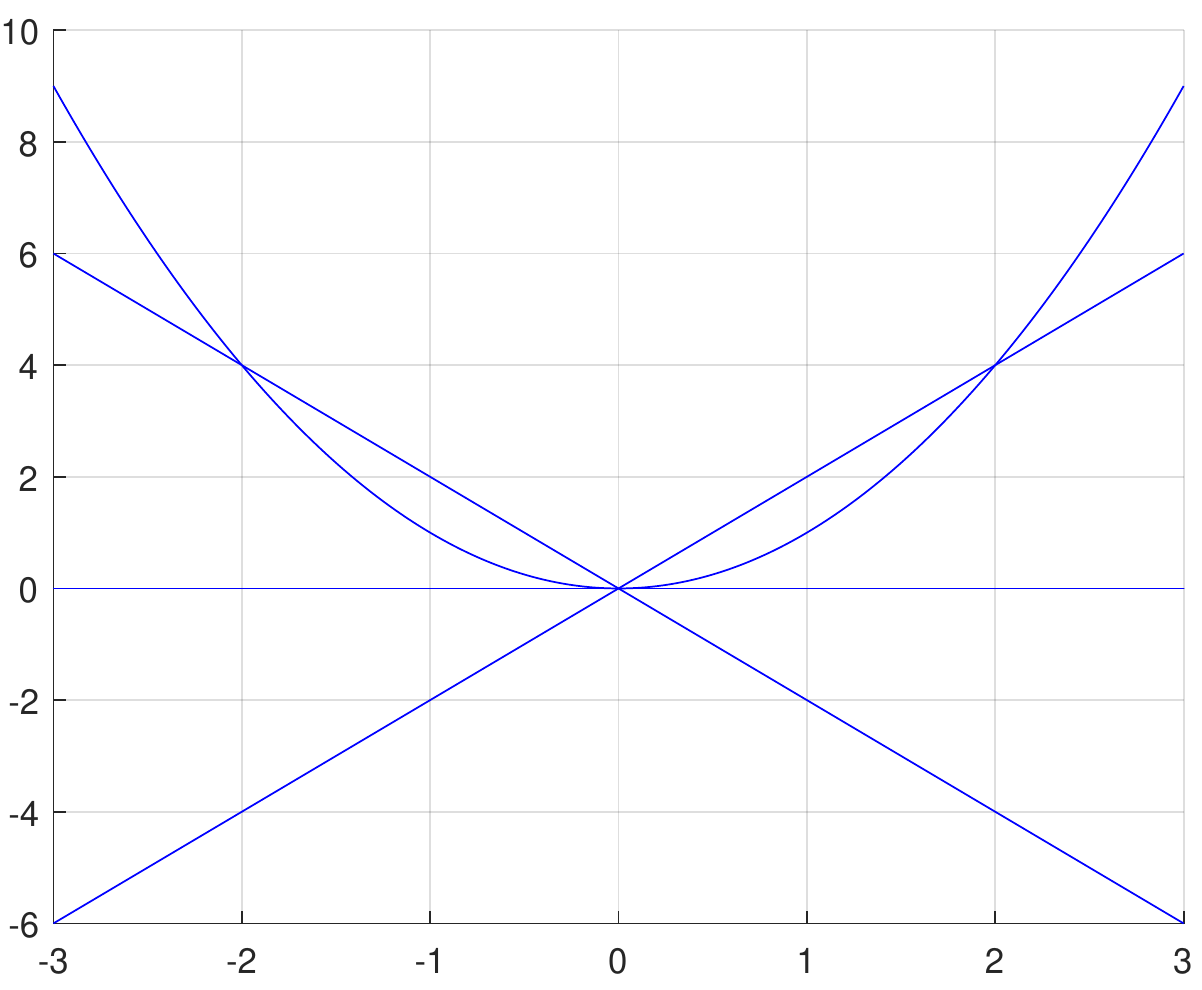}
\caption{The real picture of $(y-x^2)(y-2x)(y+2x)(y)=0$}
\label{fig:3}
\end{figure}
\end{example}

Hence, by the last example, we see that the converse of Theorem \ref{theorem} is not true in general.

Moreover, any arrangement that are constructed by a conic and two lines that are tangent to that conic has the Alexander polynomial $$\bigtriangleup(t)=(t^2+1)(t-1)^2.$$ Any such arrangement is topologically equivalent to the curve in Example \ref{example:3}. 

\begin{example}\label{example:3}
The Alexander polynomial of the curve $(y-x^2)(y+2x+1)(y-2x+1)=0$ is $(t^2+1)(t-1)^2$. Figure \ref{fig:4} is its real picture.
\begin{figure}[ht!]
\centering
\includegraphics[width=6.5cm]{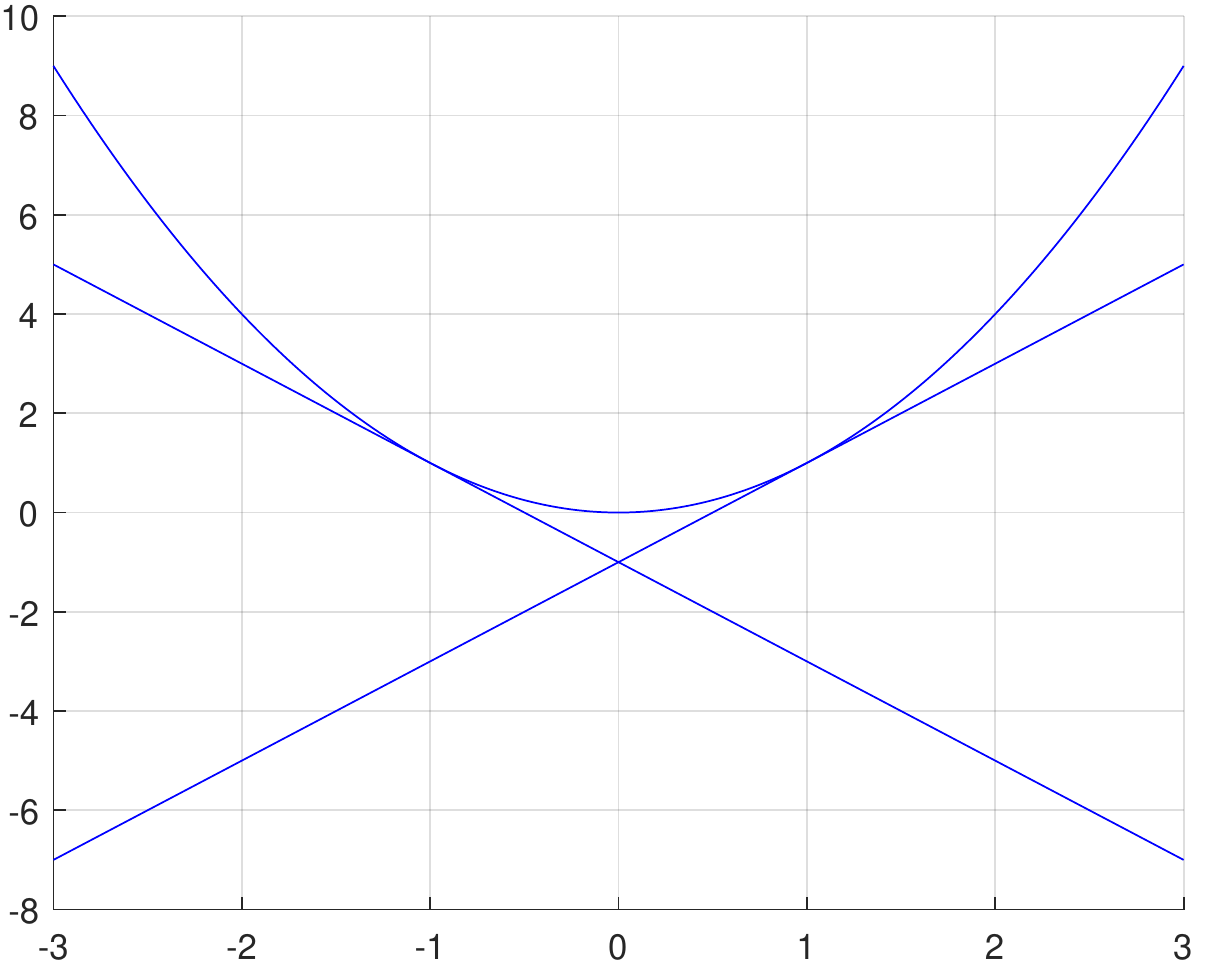}
\caption{The real picture of $(y-x^2)(y+2x+1)(y-2x+1)=0$}
\label{fig:4}
\end{figure}
\end{example}

Finally, any three curve arrangements that have only one intersection point, i.e. the curves in the form $(y+ax^n+p)(y+bx^n+p)(y+cx^n+p)=0$ with $n\in \mathbb{Z}^{+}$, $p \in \mathbb{R}$ and distinct $a,b,c\in \mathbb{R}$, also have non-trivial Alexander polynomials that is $$\bigtriangleup(t)=(t^{3n}-1)(t^{3n-3}+t^{3n-6}+...+t^3+1)(t-1).$$ See the following example.

\begin{example} 
The Alexander polynomial of the curve $(y+x^2)(y-x^2)(y)=0$ is $(t^6-1)(t^{3}+1)(t-1)$. Figure \ref{fig:5} is its real picture.
\begin{figure}[ht!]
\centering
\includegraphics[width=6.5cm]{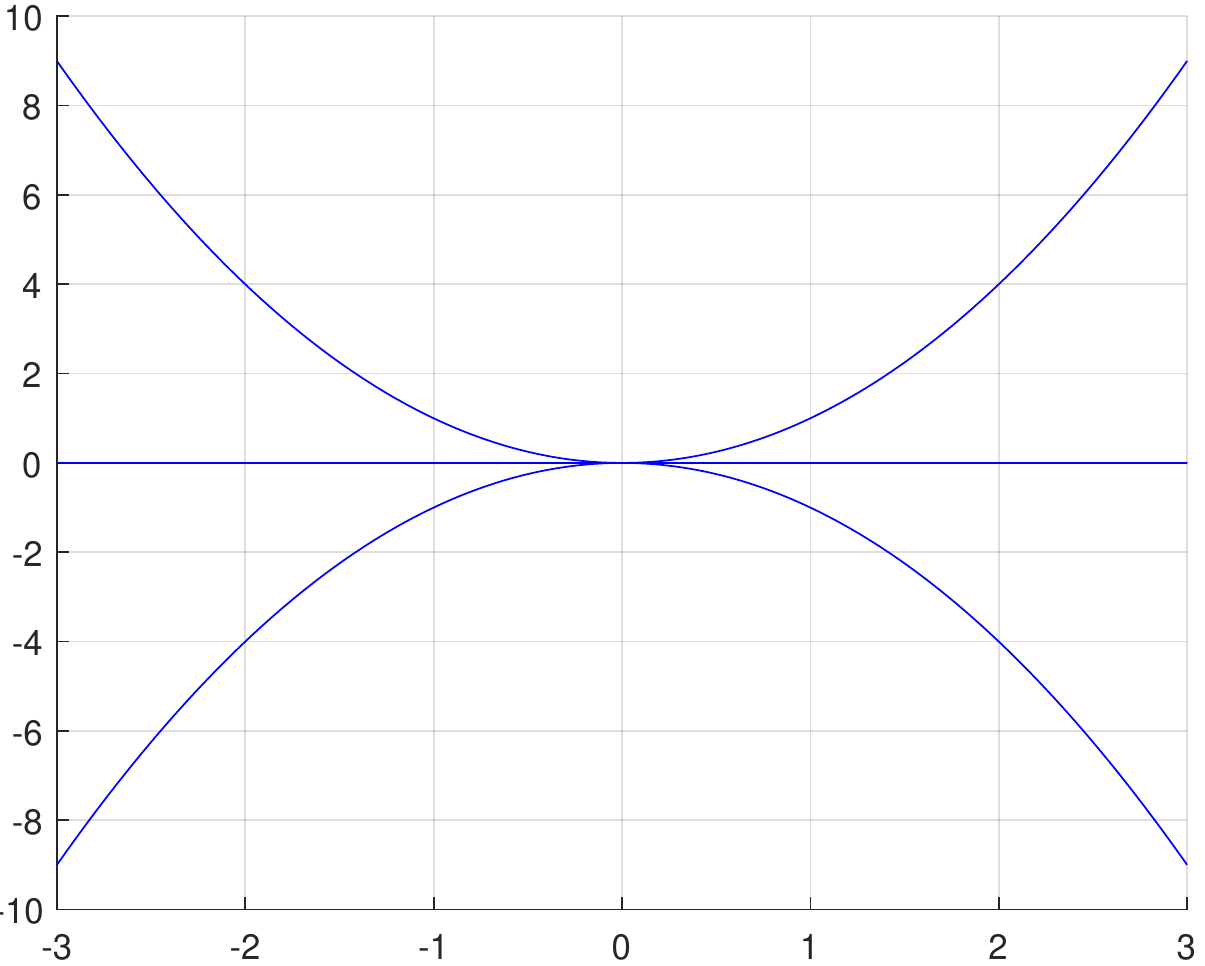}
\caption{The real picture of $(y+x^2)(y-x^2)(y)=0$}
\label{fig:5}
\end{figure}
\end{example}
\section{Comparison}
\begin{figure}[ht!]
\centering
\includegraphics[width=8cm]{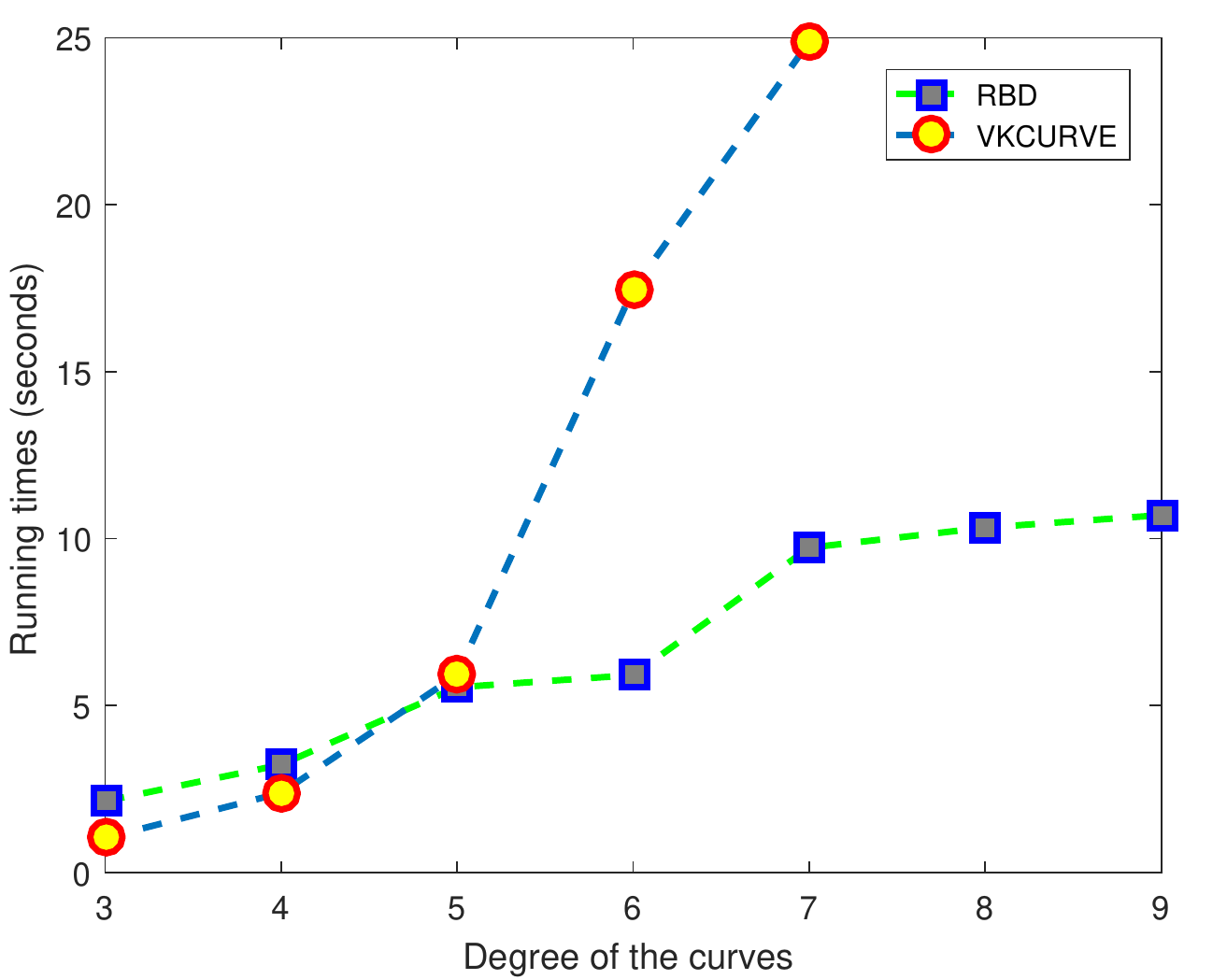}
\caption{Comparison between RBD and VKCURVE}
\label{fig:c}
\end{figure}

In this section, we compare the RBD method with VKCURVE that can be also used to compute the braid monodromy of completely reducible $n$-gonal curves. In \cite{bessis2004explicit}, the authors created the VKCURVE package which is implemented in $\mathsf{GAP}$\footnote{http://www.gap-system.org/}. We do the comparison for trigonal curves from degree 3 to 9 (we take line-line-line arrangements for degree 3, conic-line-line for degree 4 and so on). We generate random 50 polynomials for each  degree and compute the average running times for each method. Figure \ref{fig:c} shows the results.

As we see from Figure \ref{fig:c}, the gap of the run-time cost between these methods is marginal. VKCURVE could not compute the braid monodromy of the trigonal curves that have degree bigger than 7 in one day while RBD computes around 10 seconds. Furthermore, there are some trigonal curves with degree smaller than 7 where the VKCURVE could not compute the braid monodromy again in one day.

Furthermore, VKCURVE is implemented in an older version of $\mathsf{GAP}$ and cannot be used in the new version. This brings some limitations. For example, complex numbers and floating numbers are not defined in this version (only the primitive $n$-th root of unity where $n\in \mathbb{Z}$ is defined as complex numbers).

Moreover, in \cite{amoros2012monodromias}, the authors computes the braid monodromy however, this method can only compute the braid monodromy of curves up to degree 6 whereas RBD can compute for bigger degrees efficiently.

\section{Conclusion}
In this paper, we developed and implemented an efficient algorithm to compute the braid monodromy of completely reducible $n$-gonal curves. We presented another algorithm using the idea in \cite{Libgober1989} to compute the Alexander polynomials of these curves. These algorithms are implemented in Sage (available from http://www.math.fsu.edu/$\sim$maktas/ComputingBraidMonodromy). A future task is to compute other invariants based on representations of the braid group of $n$-gonal curves. 

\section*{Acknowledgments}
We would like to thank to E. Hironaka for her encouragement and many fruitful discussions on this research. We would also like to show our gratitude to Z.Y. Karatas and A. Ergur for their help with proof reading the article. We are also grateful to the anonymous reviewers for their helpful and constructive comments that greatly contributed to improving the final version of the paper. 

\providecommand{\bysame}{\leavevmode\hbox to3em{\hrulefill}\thinspace}
\providecommand{\MR}{\relax\ifhmode\unskip\space\fi MR }
\providecommand{\MRhref}[2]{%
  \href{http://www.ams.org/mathscinet-getitem?mr=#1}{#2}
}
\providecommand{\href}[2]{#2}


\begin{thebibliography}{10}

\bibitem{amoros2012monodromias}
Isabel Amor{\'o}s, Jaume andBerna~Sep{\'u}lveda, \emph{Monodrom{\'\i}as
  geom{\'e}tricas en familias de curvas de g{\'e}nero 4},  (2012).

\bibitem{amram2007fundamental}
Meirav Amram, David Garber, and Mina Teicher, \emph{Fundamental groups of
  tangent conic-line arrangements with singularities up to order 6},
  Mathematische Zeitschrift \textbf{256} (2007), no.~4, 837--870.

\bibitem{artal2001sextics}
E~Artal, J~Carmona, JI~Cogolludo, and Hiro-O Tokunaga, \emph{Sextics with
  singular points in special position}, Journal of Knot Theory and Its
  Ramifications \textbf{10} (2001), no.~04, 547--578.

\bibitem{artal2007effective}
Enrique Artal, Jorge Carmona~Ruber, and Jos{\'e} Cogolludo~Agust{\'\i}n,
  \emph{Effective invariants of braid monodromy}, Transactions of the American
  Mathematical Society \textbf{359} (2007), no.~1, 165--183.

\bibitem{artal2005topology}
Enrique Artal, Jorge~Carmona Ruber, Jos{\'e} Ignacio~Cogolludo Agust{\'\i}n,
  and Miguel~Marco Buzun{\'a}riz, \emph{Topology and combinatorics of real line
  arrangements}, Compositio Mathematica \textbf{141} (2005), no.~06,
  1578--1588.

\bibitem{bartolo2003braid}
Enrique Artal, Jorge~Carmona Ruber, and Jos{\'e}~Ignacio
  Cogolludo~Agust{\'\i}n, \emph{Braid monodromy and topology of plane curves},
  Duke Mathematical Journal \textbf{118} (2003), no.~2, 261--278.

\bibitem{artin1947theory}
Emil Artin, \emph{Theory of braids}, Annals of Mathematics (1947), 101--126.

\bibitem{arvola1992fundamental}
William~A Arvola, \emph{The fundamental group of the complement of an
  arrangement of complex hyperplanes}, Topology \textbf{31} (1992), no.~4,
  757--765.

\bibitem{bessis2004explicit}
David Bessis and Jean Michel, \emph{Explicit presentations for exceptional
  braid groups}, Experimental mathematics \textbf{13} (2004), no.~3, 257--266.

\bibitem{burau1935zopfgruppen}
Werner Burau, \emph{{\"U}ber zopfgruppen und gleichsinnig verdrillte
  verkettungen}, Abhandlungen aus dem Mathematischen Seminar der
  Universit{\"a}t Hamburg, vol.~11, Springer, 1935, pp.~179--186.

\bibitem{carmona2003monodromia}
J~Carmona, \emph{Monodrom{\i}a de trenzas de curvas algebraicas planas},
  Thesis, Universidad de Zaragoza (2003).

\bibitem{chisini1933suggestiva}
Oscar Chisini, \emph{Una suggestiva rappresentazione reale per le curve
  algebriche piane}, Ist. Lombardo, Rend., II. Ser \textbf{66} (1933),
  1141--1155.

\bibitem{chisini1952courbes}
\bysame, \emph{Courbes de diramation des plans multiples et tresses
  alg{\'e}briques}, Deuxieme Colloque de G{\'e}om{\'e}trie Alg{\'e}brique,
  Masson, Paris, 1952, pp.~11--27.

\bibitem{cordovil1995braid}
R~Cordovil and JL~Fachada, \emph{Braid monodromy groups of wiring diagrams},
  Bollettino della Unione Matematica Italiana-B (1995), no.~2, 399--416.

\bibitem{deconinck2001computing}
Bernard Deconinck and Mark Van~Hoeij, \emph{Computing \text{R}iemann matrices
  of algebraic curves}, Physica D: Nonlinear Phenomena \textbf{152} (2001),
  28--46.

\bibitem{degtyarev2012topology}
Alex Degtyarev, \emph{Topology of algebraic curves: an approach via dessins
  d'enfants}, vol.~44, Walter de Gruyter, 2012.

\bibitem{hironaka1993abelian}
Eriko Hironaka, \emph{Abelian coverings of the complex projective plane
  branched along configurations of real lines}, no. 502, American Mathematical
  Soc., 1993.

\bibitem{eko}
\bysame, \emph{{Computing Alexander polynomials using monodromy}},  (2011),
  unpublished.

\bibitem{kaplan2004braid}
Shmuel Kaplan, Eran Liberman, and Mina Teicher, \emph{Braid monodromy
  computation of real singular curves}, arXiv preprint math/0410444 (2004).

\bibitem{kharlamov2003braid}
Viatcheslav~Mikhailovich Kharlamov and Vik~S Kulikov, \emph{On braid monodromy
  factorizations}, Izvestiya: Mathematics \textbf{67} (2003), no.~3, 499.

\bibitem{kulikov2000braid}
Vik~S Kulikov and Mina Teicher, \emph{Braid monodromy factorizations and
  diffeomorphism types}, Izvestiya: Mathematics \textbf{64} (2000), no.~2, 311.

\bibitem{libgober1986homotopy}
Anatoly Libgober, \emph{On the homotopy type of the complement to plane
  algebraic curves}, J. reine angew. Math \textbf{367} (1986), 103--114.

\bibitem{Libgober1989}
\bysame, \emph{{Invariants of plane algebraic curves via representations of the
  braid groups}}, Inventiones mathematicae \textbf{30} (1989), no.~6, 25--30.

\bibitem{moishezon1981stable}
Boris Moishezon, \emph{Stable branch curves and braid monodromies}, Algebraic
  geometry, Springer, 1981, pp.~107--192.

\bibitem{moishezon1991braid}
Boris Moishezon and Mina Teicher, \emph{Braid group technique m complex
  geometry, ii: From arrangements of lines and conics to cuspidal curves},
  Algebraic geometry, Springer, 1991, pp.~131--180.

\bibitem{oka1992symmetric}
Mutsuo Oka, \emph{Symmetric plane curves with nodes and cusps}, Journal of the
  Mathematical Society of Japan \textbf{44} (1992), no.~3, 375--414.

\bibitem{oka1996two}
\bysame, \emph{Two transforms of plane curves and their fundamental groups},
  Journal of Mathematical Sciences-University of Tokyo \textbf{3} (1996),
  no.~2, 399--444.

\bibitem{Oka2005}
\bysame, \emph{{A survey on Alexander polynomials of plane curves}},
  Singularit\'{e}s Franco-Japonaise, S\'{e}minaire et congres (2005), 209--232.

\bibitem{oka2002fundamental}
Mutsuo Oka and Duc~Tai Pho, \emph{Fundamental group of sextics of torus type},
  Trends in singularities, Springer, 2002, pp.~151--180.

\bibitem{salvetti1987topology}
Mario Salvetti, \emph{Topology of the complement of real hyperplanes in ℂ n},
  Inventiones mathematicae \textbf{88} (1987), no.~3, 603--618.

\bibitem{csimcsek2004computer}
Hakan {\c{S}}im{\c{s}}ek, Mustafa Bayram, and Ugur Yavuz, \emph{A computer
  program to calculate alexander polynomial from braids presentation of the
  given knot}, Applied mathematics and computation \textbf{153} (2004), no.~1,
  199--204.

\bibitem{van1933fundamental}
Egbert~R Van~Kampen, \emph{On the fundamental group of an algebraic curve},
  American journal of Mathematics (1933), 255--260.

\bibitem{zariski1929problem}
Oscar Zariski, \emph{On the problem of existence of algebraic functions of two
  variables possessing a given branch curve}, American Journal of Mathematics
  (1929), 305--328.

\end{thebibliography}
\end{document}